\documentclass[12pt,a4paper]{amsart}
\usepackage{amsmath, amssymb, amsthm,
amsfonts, amsthm, bm, eucal, graphicx}
\usepackage[all]{xy}
\usepackage{graphics}
\usepackage{epsfig}
\usepackage{graphicx}
\theoremstyle{plain}
\usepackage{amssymb}
\usepackage{enumerate}

\advance\hoffset-10mm \advance\textwidth40mm

\newtheorem{thm}{Theorem}[section]
\newtheorem{cor}[thm]{Corollary}
\newtheorem{lem}[thm]{Lemma}
\newtheorem{prop}[thm]{Proposition}

\theoremstyle{definition}
\newtheorem{defi}[thm]{Definition}
\newtheorem{defis}[thm]{Definitions}
\newtheorem{conj}[thm]{Problem}
\newtheorem{conv}[thm]{Convention}
\newtheorem{nota}[thm]{Notation}
\newtheorem{rem}[thm]{Remark}
\newtheorem{rems}[thm]{Remarks}
\newtheorem{exa}[thm]{Example}
\newtheorem{exas}[thm]{Examples}
\newtheorem{sit}[thm]{}

\newcommand{\brem}{\begin{rem}}
\newcommand{\brems}{\begin{rems}}
\newcommand{\erem}{\end{rem}}
\newcommand{\erems}{\end{rems}}
\newcommand{\bexa}{\begin{exa}}
\newcommand{\bexas}{\begin{exas}}
\newcommand{\eexa}{\end{exa}}
\newcommand{\eexas}{\end{exas}}
\newcommand{\bdefi}{\begin{defi}}
\newcommand{\edefi}{\end{defi}}
\newcommand{\bdefis}{\begin{defis}}
\newcommand{\edefis}{\end{defis}}
\newcommand{\bcor}{\begin{cor}}
\newcommand{\ecor}{\end{cor}}
\newcommand{\blem}{\begin{lem}}
\newcommand{\elem}{\end{lem}}
\newcommand{\bconv}{\begin{conv}}
\newcommand{\econv}{\end{conv}}
\newcommand{\bconj}{\begin{conj}}
\newcommand{\econj}{\end{conj}}
\newcommand{\bprop}{\begin{prop}}
\newcommand{\eprop}{\end{prop}}
\newcommand{\bthm}{\begin{thm}}
\newcommand{\ethm}{\end{thm}}
\newcommand{\bnota}{\begin{nota}}
\newcommand{\enota}{\end{nota}}
\newcommand{\bsit}{\begin{sit}}
\newcommand{\esit}{\end{sit}}
\newcommand{\be}{\begin{equation}}
\newcommand{\ee}{\end{equation}}
\newcommand{\bproof}{\begin{proof}}
\newcommand{\eproof}{\end{proof}}
\def\ba{\begin{array}}
\def\ea{\end{array}}
\def\bea{\begin{eqnarray}}
\def\eea{\end{eqnarray}}

\newtheorem*{theo*}{Theorem}

\theoremstyle{definition}

\newtheorem*{definition*}{Definition}

\theoremstyle{definition}
\newtheorem{case}[thm]{}
\newcommand\bys[1]{\begin{case}{\bf #1} }
\newcommand\eys{\end{case}}

\newcommand{\LND}{\operatorname{LND}}
\newcommand{\Der}{\operatorname{Der}}

\newcommand{\Proj}{\operatorname{Proj}}
\newcommand{\Frac}{\operatorname{Frac}}
\newcommand{\spec}{\operatorname{Spec}}
\newcommand{\supp}{\operatorname{supp}}

\newcommand{\affcone}{\operatorname{cone}}
\newcommand{\Pic}{\operatorname{Pic}}

\newcommand{\Div}{\operatorname{Div}}

\newcommand{\Sing}{\operatorname{Sing}}

\newcommand{\im}{\operatorname{im}}

\def\fp{{\mathfrak p}}

\def\cO{{\mathcal O}}

\def\p{{\partial}}

\newcommand{\A}{{\mathbb A}}
\newcommand{\PP}{{\mathbb P}}
\newcommand{\DD}{{\mathbb D}}

\newcommand{\C}{{\mathbb C}}
\newcommand{\Q}{{\mathbb Q}}
\newcommand{\Z}{{\mathbb Z}}
\newcommand{\N}{{\mathbb N}}

\newcommand{\V}{{\mathbb V}}

\newcommand{\G}{{\mathbb G}}
\newcommand{\kk}{{\Bbbk}}

\def\cO{{\mathcal O}}

\def\bdi{\pdfsyncstop\begin{diagram}}
\def\edi{\end{diagram}\pdfsyncstart}

\title[$\G_{\mathrm a} $-actions on affine cones]{$\G_{\mathrm a} $-actions
on affine cones}

\author{Takashi Kishimoto}

\address{Department of Mathematics,
Faculty of Science, Saitama University, Saitama 338-8570, Japan}
\email{tkishimo@rimath.saitama-u.ac.jp}

\author{Yuri Prokhorov}
\thanks{The first author was supported by a Grant-in-Aid
for Scientific Research of JSPS, No. 24740003.
The second author was  partially supported by
RFBR grant No. 11-01-00336-a, the grant of
Leading Scientific Schools No. 4713.2010.1,
Simons-IUM fellowship, and
AG Laboratory SU-HSE, RF government
grant ag. 11.G34.31.0023. 
}
\address{Steklov Institute of Mathematics,
8 Gubkina street, Moscow 119991
and
Laboratory of Algebraic Geometry, GU-HSE, 7 Vavilova street,
Moscow 117312}
\email{prokhoro@gmail.com}

\author{Mikhail Zaidenberg}

\address{Universit\'e
Grenoble I, Institut Fourier, UMR 5582 CNRS-UJF, BP 74, 38402 St.\
Martin d'H\`eres c\'edex, France} \email{zaidenbe@ujf-grenoble.fr}

\thanks{}

\begin{document}

\begin{abstract}
An affine algebraic variety $X$ is called {\em cylindrical} if it
contains a principal Zariski dense open cylinder $U\simeq
Z\times\A^1$.  A polarized projective variety $(Y,H)$  is called
{\em cylindrical} if it contains a cylinder $U=Y\setminus \supp
D$, where $D$ is an effective $\Q$-divisor on $Y$ such that  $[D]\in \Q_{+} [H]$
in $\Pic_\Q(Y)$. We
show that cylindricity of a polarized projective variety
is equivalent to that of a certain Veronese affine cone over this variety.
This gives a criterion of existence of a unipotent group action on an affine
cone. 
\end{abstract}

\maketitle

\bigskip 


\section*{Introduction} Fix an algebraically closed field $\kk$ of
characteristic zero, and let $\G_{\mathrm a}=\G_{\mathrm a}
(\kk)$, the additive group of $\kk$. We investigate when  the
affine cone over an irreducible, normal projective variety over
$\kk$ admits a non-trivial action of a unipotent group. Since any
unipotent group contains a one parameter unipotent subgroup,
instead of considering general unipotent group actions we stick to
the $\G_{\mathrm a} $-actions. Our main purpose in this paper is
to provide a geometric criterion for existence of such an action,
see Theorem \ref{ext-grad-int} and Corollary \ref{new-cor}. The
former version of such a criterion
 in \cite{KPZ1} involved some unnecessary
assumptions. In Theorem \ref{ext-grad-int} we remove these
 assumptions. What is more important, we extend  our  criterion
 so that it can be applied more generally to affine
 quasicones. An affine
 quasicone is an affine variety $V$ equipped with a $\G_{\mathrm m} $-action
such that the fixed point set $V^{\G_{\mathrm m} }$ attracts the whole $V$.
Thus the variety  $Y=(V\setminus V^{\G_{\mathrm m} })/\G_{\mathrm m} $
is projective over the
affine variety $S=V^{\G_{\mathrm m} }$. We assume in this paper that $Y$ is
normal. Our criterion is formulated in terms of a geometric
property called {\em cylindricity}, which merits study for its own sake.

\bys{Cylindricity.}
Let us fix the notation. For two $\Q$-divisors $H$ and $H'$
on a quasiprojective variety
$Y$ we write
$H\sim H'$ if $H$ and $H'$ are linearly equivalent that is,
$H-H'={\rm div}\, (f)$
for a rational function $f$ on $Y$.
We write $[H']\in \Q_{+} [H]$ in $\Pic_\Q(Y)$ meaning that
$H'\sim\frac{p}{q}H$ for
some coprime positive integers $p$ and $q$.

\begin{defi}[cf.\ {\cite[3.1.4]{KPZ1}}]
\label{cyl} Let $Y$ be  a quasiprojective variety over $\kk$
polarized by an  ample $\Q$-divisor $H\in \Div_\Q(Y)$. We say that
the pair $(Y, H)$ is {\em cylindrical} if there exists an
effective $\Q$-divisor $D$ on $Y$ such that $[D]\in \Q_{+} [H]$ in
$\Pic_\Q(Y)$ and $U=Y\setminus\supp D$ is a {\em cylinder} i.e.
$$ U\simeq Z\times\A^1\,$$ for some variety $Z$.
Here $U$ and $Z$ are quasiaffine varieties. Such a cylinder $U$ is
called {\em $H$-polar} in \cite[3.1.7]{KPZ1}.
Notice that the cylindricity of $(Y, H)$
depends only on
the ray $\Q_+[H]$ generated by $H$ in $\Pic_\Q Y$.
\end{defi}

\brem\label{0001}
The pair $(Y, H)$ can admit several essentially different  cylinders.
For instance, let $Y=\PP^1$
and $H$ is a $\Q$-divisor on $Y$ of positive degree.
Then any  divisor $D=rP$, where $P\in\PP^1$ and $r\in\Q_+$,
defines an $H$-polar
cylinder on $Y$.

\erem

\bdefi\label{cyl-aff} An affine variety $X$ is called  {\em
cylindrical} if it contains a principal cylinder
$$\DD(h):=X\setminus \V(h)\simeq Z\times\A^1,\quad\mbox{where}\quad
\V(h)=h^{-1}(0),\,$$ for some variety $Z$ and some
regular function $h\in\cO(X)$. Hence $U$ and $Z$ are  affine varieties.
 \edefi

The  cylindricity of affine varieties is important due to the
following well known fact (see e.g.\ \cite[Proposition
3.1.5]{KPZ1}).

\bprop\label{ext} An
affine variety  $X=\spec A$ over $\kk$  is
cylindrical  if and only if it admits
an effective $\G_{\mathrm a} $-action, if and only if $\LND(A)\neq\{0\}$,
where $\LND(A)$
stands for the set of all 
locally
nilpotent derivations on $A$.
\eprop

The proof is based upon the slice construction, which we recall in
subsection \ref{Slice-construction}. In Section
\ref{section-Preliminaries} we gather necessary preliminaries on
positively graded rings. In particular, we   give a graded version
of the  slice construction, and recall the DPD
(Dolgachev-Pinkham-Demazure) presentation of a positively graded
affine domain $A$ over $\kk$ in terms of an ample
$\Q$-divisor $H$ on the variety $Y=\Proj A$. In Section
\ref{section-The-criterion} we prove  our main result,
inspired by Theorem 3.1.9 in \cite{KPZ1}.

\bthm\label{ext-grad-int} 
Let $A=\bigoplus_{\nu\ge 0} A_\nu$
be a positively graded affine domain over $\kk$.
Define the projective  variety $Y=\Proj A$
relative to this grading, let $H$  be the associated
$\Q$-divisor on $Y$, and let $V=\spec A$,
the affine quasicone over $Y$.
\begin{enumerate}\item[(a)] If $V$ is
cylindrical, then the associated pair $(Y,H)$ is  cylindrical.
\item[(b)] If the pair
$(Y,H)$ is cylindrical, then for some $d\in\N$ the Veronese cone
$V^{(d)}=\spec A^{(d)}$ is cylindrical, where
$A^{(d)}=\bigoplus_{\nu\ge 0} A_{d\nu}$. 
\end{enumerate}
\ethm

In Lemma \ref{very-ample}, we specify the range of values of
$d$ satisfying the second assertion. In particular, the
latter holds with $d=1$
provided that $H\in\Div (Y)$. Thus if $H\in\Div (Y)$,
Theorem  \ref{ext-grad-int} yields a necessary and sufficient condition
of cylindricity of the corresponding affine quasicone,
see Corollary \ref{new-cor}.

 Let us notice that the main difference between Theorem
\ref{ext-grad-int} and the former criterion in \cite[Theorem
3.1.9]{KPZ1} consists in removing the unnecessary extra assumption
$\Pic(Z)=0$ or, what is equivalent, $\Pic(U)=0$ used in the proof
in \cite{KPZ1}. So our proof of Theorem \ref{ext-grad-int} here is
pretty much different. On the other hand, it is worthwhile
mentioning that the present proof works only under the assumption
of normality of the variety $Y$ (see \S \ref{ample} below), absent
in \cite[Theorem 3.1.9]{KPZ1} and as well in Proposition
\ref{ext}. 

In Section \ref{section-Final-remarks}
we provide several examples that illustrate our criterion.
Besides, we discuss a possibility to lift a $\G_{\mathrm a} $-action
on a Veronese cone $V^{(d)}$ over $Y$ to the affine cone $V$ over $Y$. 
In particular, we prove the criterion of Corollary \ref{new-cor}
cited above. 

It is our pleasure  to thank Michel Brion, Hubert Flenner,
Shulim Kaliman, Kevin Langlois,  Alvaro Liendo, and Alexander
Perepechko for useful discussions and references. The discussions
with Alvaro Liendo and Alexander Perepechko were especially
pertinent and allowed us to improve significantly the
presentation. Our thanks are also due to the referee for his
very helpful remarks and comments. After our paper \cite{KPZ4} has been published, 
Kevin Langlois pointed out that the argument in the proof of \cite[Lemma \ref{spec-case}]{KPZ4} needs to be corrected. We are grateful to Kevin Langlois for this observation.
Certainly, this lemma is not essential in the proof of our main 
Theorem \ref{ext-grad-int}. 
Indeed, this proof exploits Proposition \ref{propo}, which does not depend on Lemma \ref{spec-case} and deals with 
the general case in Theorem \ref{ext-grad-int}(a). However, the particular case of this theorem considered in Lemma \ref{spec-case} and, especially, in
Corollary \ref{cyl-fac} is the most useful in applications. Hence we provide here 
a corrected proof of Lemma \ref{spec-case}.

\eys

\section{Preliminaries}\label{section-Preliminaries}
Throughout this article, $A$ will denote an affine domain over $\kk$, and $\LND(A)$ will denote the set of locally nilpotent derivations of $A$.

\bys{Slice construction.}\label{Slice-construction}
Let $\p\in \LND(A)$ be non-zero.
The filtration \be\label{gsc} A^\p=\ker\p\varsubsetneq \ker
\p^2\varsubsetneq\ker \p^{3}\varsubsetneq \ldots\ee being strictly increasing
one can find an element $g\in \ker \p^{2}\setminus\ker \p$.
Letting $h=\p g\in \ker\p\cap {\rm im}\,\p$, where $h\neq 0$, one
considers the localization $A_{h}=A[h^{-1}]$ and the principal
Zariski dense open subset
$$\DD(h)=X\setminus\V(h)\simeq\spec
A_{h}\,,\quad\mbox{where}\quad \V(h)=h^{-1}(0)\,.$$
The derivation $\p$ extends to a locally nilpotent
derivation on $A_{h}$ denoted by the same letter. The element
$s=g/h\in A_{h}$ is a slice of $\p$ that is, $\p (s)=1$. Hence
$$A_{h}=A_{h}^\p[s],\qquad\mbox{where}\quad\p=d/d s\quad\mbox{and}
\quad A_{h}^\p\simeq A_{h}/(s)\,$$ (`Slice Theorem', \cite[Corollary
1.22]{Fr}). Thus $ \DD(h)\simeq Z\times\A^1$ is a principal cylinder
in $X$ over $Z=\spec A_{h}^\p$. The $\G_{\mathrm a} $-action on $\DD(h)$
associated with $\p$ is defined by the translations along the
second factor. The natural projection $p_1:\DD(h)\to Z$ identifies
$\V(g)\setminus\V(h)\subseteq \DD(h)$ with $Z$. Choosing $f\in
A_{h}^\p=\cO(Z)$ such that $\Sing(Z)\subseteq\V(f)$ we can replace
$g$ and $h$ by $fg$ and $fh$, respectively, so that the slice $s$
remains the same, but the new cylinder $\DD(fh)$ over an affine
variety $Z'=\DD(fh)/\G_{\mathrm a} $ is smooth.

\bsit\label{homog} {\bf Graded slice construction.}
Suppose that the ring $A$ is graded. 
The gradings used in this paper are $\Z$-gradings
over $\kk$, i.e., if $A=\bigoplus_{\nu\in\Z} A_{\nu}$
is a $\Z$-grading of the ring $A$,
then $\kk \subseteq A_0$. The grading is said to be a
{\em positive grading} if $A_\nu=\{0\}$
for $\nu< 0$; we do not assume that $A_0=\kk$. 

 Any non-zero derivation $\eta\in\LND(A)$ can be decomposed 
into a sum of homogeneous components
$$\eta=\sum_{i=1}^n \eta_i,\quad\mbox{where}
\quad\eta_i\in\Der(A),\quad\deg \eta_i<\deg \eta_{i+1}\quad\forall
i,\quad\mbox{and}\quad \eta_n\neq 0\,.$$ Letting $\p=\eta_n$
be the principal homogeneous component of $\eta$, 
$\p$ is again
locally nilpotent and homogeneous (see \cite{Da, Re}). Hence all
kernels in (\ref{gsc}) are graded.
 So one can choose homogeneous
elements $g$, $h$, and $s$ as in \S \ref{Slice-construction}.
With this choice, we call the
construction of the cylinder in \S \ref{Slice-construction}
a {\em graded slice construction}. 
\esit
\eys

\bys{Graded rings and associated
schemes.}\label{gras}
We recall some well known facts on positively graded rings and
associated schemes. The presentation
below is borrowed from  \cite{De}, \cite[sect.\ 2]{Fl}, \cite[\S 2.1]{FZ1},
and
\cite[Lecture 3]{Do}.

\bnota\label{a.1} Given a graded affine domain
$A=\bigoplus_{\nu\in\Z} A_{\nu}$ over $\kk$ the group $\G_{\mathrm
m} $ acts on $A$ via $t.a=t^{\nu}a$ for $a\in A_{\nu}$. This
action is effective if and only if the {\em saturation index}
$e(A)$ equals 1, where $$e(A)=\gcd\{\nu\,|\,A_\nu\neq (0)\}\,.$$
If $A$ is positively  graded 
 then the associated
scheme $Y=\Proj A$
is projective  over\footnote{Notice
that $A_0$ can be here an arbitrary affine domain.}
the affine scheme $S=\spec
A_0$ \cite[II]{EGA}. Furthermore, $Y$ is covered by the affine
open subsets 
$$
\DD_+(f)=\{\fp\in \Proj A : f\notin \fp\}\cong \spec
A_{(f)}\,,
$$ 
where $f\in A_{>0}$ is a homogeneous element and
$A_{(f)}=(A_{f})_0$ stands for the degree zero part of the
localization $A_f$. The affine variety $V=\spec A$ is called a
{\em quasicone} over $Y$ with {\em vertex} $\V(A_{>0})$ and
with {\em punctured quasicone} $V^*=V \backslash \V(A_{>0})$,
where $\V(I)$ stands for the zero set of an ideal $I\subseteq A$.
For a homogeneous ideal $I\subseteq A$, $\V_+(I)$ stands for its
zero set in $Y=\Proj A$. There is a natural surjective morphism
$\pi:V^*\to Y$. If $A=A_0[A_1]$ i.e., $A$ is generated as an
$A_0$-algebra by the elements of degree $1$ then $V^*\to Y=\Proj A$
is a locally trivial $\G_{\mathrm m} $-bundle. In the general case the
following holds. \enota

\blem\label{sl} $\Proj  A\cong V^*/\G_{\mathrm m} $. \elem

\bproof Indeed, $V^*$ is covered by the $\G_{\mathrm m} $-invariant affine
open subsets $\DD(f)=\spec A_f$, where $f\in A_{d}$ with $d >0$.
Since $(A_f)^{\G_{\mathrm m} }=A_{(f)}=(A_f)_0$ we have
 $\DD_+(f)=\DD(f)/{\G_{\mathrm m} }$ and the lemma follows. \eproof

\brem\label{trivialization} Assuming that $A$ is a domain over $\kk$ and
$e(A)=1$ one can find a pair of non-zero homogeneous elements $a\in
A_\nu$ and $b\in A_\mu$ of coprime degrees. Let $p, q\in \Z$ be
such that $p\nu+q\mu=1$. Then the localization $A_{ab}$ is graded,
the element $u=a^pb^q\in (A_{ab})_1$ is invertible, and
$A_{ab}=A_{(ab)}[u,u^{-1}]\,.$ This gives a trivialization of the
orbit map $\pi:V^*\to  Y=\Proj  A$ over the principal open set
$\DD_+(ab)\subseteq Y$: $$\DD(ab)=\pi^{-1}(\DD_+(ab))\simeq
\DD_+(ab)\times \A^1_*,\quad\mbox{where}\quad
\A^1_*=\A^1\setminus\{0\}\,.$$ \erem

\bsit\label{fle}  {\bf Cyclic quotient construction.}
Let $h\in A_{m}$ be a homogeneous element of degree $m>0$,
and let $F=A/(h-1)$.
For  $a\in A$ we let $\bar a$ denote the class of $a$ in $F$.
The  projection $\rho:A\to F$, $a\mapsto \bar a$,
extends to the localization $A_h$ via
$\rho(a/h^l)=\rho(a)=\bar a$.
The cyclic group
$\boldsymbol\mu_{m}\subseteq\G_{\mathrm m} $ of the
$m$th roots of unity acts on $F$ effectively and so defines a $\Z_m$-grading
$$F=\bigoplus_{[i]\in\Z_m} F_{[i]}\,,$$
where $\Z_m=\Z/m\Z$ and $[i]\in\Z_m$ stands for the residue class of
$i\in\Z$ modulo $m$.
It is easily seen that the morphism $\rho: A_{h}\to F$ restricts to
an isomorphism
$\rho: A_{(h)}\stackrel{\simeq}{\longrightarrow} F_{[0]}$.
This yields  a cyclic quotient
$$Y_h\to Y_h/\boldsymbol\mu_m=\DD_+(h)\subseteq Y\,,\,\,\,\,\mbox{where}\,\,\,\,
Y_h=h^{-1}(1)\subseteq V\,,$$
$$V=\spec A,\quad\mbox{and}\quad\DD_+(h)=\spec A_{(h)}\simeq \spec F_{[0]}\,.$$

Let  $\p$ be a homogeneous locally nilpotent derivation of $A$. If $h\in A^\p_m$ then
the principal ideal $(h-1)$ of $A$ is $\p$-stable. Hence the hypersurface $Y_h=\V(h-1)$
is stable
under the $\G_{\mathrm a} $-action on $V$
generated by $\p$,
and $\p$ induces
a homogeneous locally nilpotent derivation $\bar\p$ of the $\Z_m$-graded ring $F$.
The kernel $F^{\bar\p}=\ker\bar\p$ is a  $\Z_m$-graded subring of $F$:
$$F^{\bar \p}=\bigoplus_{[i]\in\Z_m} F^{\bar \p}_{[i]}\,,
\quad\mbox{where}\quad F^{\bar \p}_{[i]}=F_{[i]}\cap F^{\bar \p}\,.$$

Assume further  that $F$ is a domain.
Then the set $\{[i]\in\Z_m\,|\,F_{[i]}\neq (0)\}$
is a cyclic subgroup, say,
$\boldsymbol\mu_n\subseteq\boldsymbol\mu_m$.
Letting $k=m/n$ we can write
$$F^{\bar \p}=\bigoplus_{i=0}^{n-1} F^{\bar \p}_{[ki]}\,.$$
\esit

\blem\label{lie} We have
$k=e(A^{\p}) (=e(A_h^\p))$ and
$\gcd (k,d)=1$, where $d=-\deg\p$. \elem

\bproof The second assertion follows from the first since $\gcd(d,\,e(A^\p))=1$.
Indeed,  notice that for any non-zero homogeneous element 
$g\in A_j$ ($j>0$) there is $r\in\N$ such that 
$\p^{r}g\in A^\p_{j-rd}\setminus\{0\}$. Hence $j=rd+se(A^\p)$ for some $s\in\Z$.
Since by our assumption $e(A)=1$ it follows that $\Z=\langle d,\,e(A^\p)\rangle$
and so  $\gcd(d,\,e(A^\p))=1$.

To prove the first equality we let $g\in A_j$ be such that $\bar
g\in F^{\bar \p}$, where $j>0$. The restriction $g|_{Y_h}$ being
invariant under the induced $\G_{\mathrm a} $-action on $Y_h$,
this restriction is constant on any $\G_{\mathrm a} $-orbit in
$Y_h$. For a general point $x \in \DD (h) \subseteq V$ there
exists $\lambda\in \G_{\mathrm m} $ such that $h(\lambda.x)=1$.
Since $\p$ is homogeneous the $\G_{\mathrm m} $-action on $V$
induced by the grading normalizes  the $\G_{\mathrm a} $-action.
Therefore $\lambda.(\G_{\mathrm a} .x)=\G_{\mathrm a}
.(\lambda.x)\subseteq Y_h$ and so $g|_{\G_a .(\lambda.x)}$ 
is constant. Hence also $g|_{\lambda.(\G_{\mathrm a}
.x)}=\lambda^j g|_{\G_{\mathrm a} .x}$ is. It follows that $g\in
A_j^\p$. Clearly $\rho (A^\p)\subseteq F^{\bar\p}$, so finally
$\rho (A^\p)= F^{\bar\p}$. Thus $k=e(A^{\p}) $.
 \eproof

\eys

\bys{Quasicones and ample $\Q$-divisors.}\label{ample} To any pair
$(Y, H)$, where  $Y\to S$ is a proper  normal integral
$S$-scheme, $S=\spec A_0$ is a normal affine variety
over $\kk$, and $H$ is an ample $\Q$-divisor on $Y$, one can
associate a positively graded integral domain over $\kk$,
\be\label{pogr} A=A(Y,H)=\bigoplus_{\nu\ge 0}
A_\nu,\quad\mbox{where}\quad A_\nu= H^0(Y, \cO_Y(\lfloor\nu
H\rfloor))\,.\ee The algebra $A$ has saturation index $e(A)=1$, is
finitely generated and normal\footnote{See \cite[Ch.\ II, Exercise
5.14(a)]{Ha}, \cite[3.1]{De}, \cite[3.3.5]{Do}, and \cite[Theorem
3.1]{AH}.}. So the associated affine quasicone $V=\spec A$ over
$Y=\Proj A$ is normal.

 Conversely, every affine quasicone $V=\spec A$, where $A$ is a
normal affine positively graded $\kk$-domain of
dimension at least 2 and with  saturation index $1$
 arises in this way (\cite[3.5]{De}). The corresponding ample $\Q$-divisor $H$
on $Y$ is defined uniquely by the quasicone $V$ up
to the linear equivalence\footnote{See \cite[Theorem 3.3.4]{Do}; cf.\ also
\cite{FZ3} and
\cite[Theorem 3.4]{AH}.}.
In particular, its fractional part $\{H\}=H-\lfloor H\rfloor$
is uniquely determined by $V$; it is called
the {\em Seifert divisor} of the quasicone $V$, see \cite[3.3.2]{Do}.

\bsit\label{triv}  Let again $A=A(Y,H)$ be as in (\ref{pogr}).
By virtue of Remark \ref{trivialization} there exists  on $V$
a homogeneous rational function $u\in ({\rm Frac}\, A)_1$ of degree 1.
Notice that the divisor ${\rm div}\,u$ on $V$ is $\G_{\mathrm m} $-invariant.
Choosing this function suitably one can achieve that
${\rm div} (u|_{V^*})=\pi^*(H)$, where $\pi:V^*\to Y$
is the quotient by the $\G_{\mathrm m} $-action (see Lemma \ref{sl}).

Furthermore,
$ {\rm Frac}\, A= ({\rm Frac}\, A)_0(u)$. So
any homogeneous   rational function $f\in (\Frac A)_d$
of degree $d$ on $V$ can be written as $\psi u^d$ for some
$\psi \in ({\rm Frac}\, A)_0$.
For any   $d>0$ the $\Q$-divisor class $[dH]$ is ample.
It  is invertible
and trivial on any open set $\DD_+(ab)\subseteq Y$ as in
Remark \ref{trivialization}.

A rational function on $Y$ can be lifted to a
$\G_{\mathrm m} $-invariant rational function on $V$.
Thus the field $({\rm Frac}\, A)_0$ can be naturally identified
with the function field of $Y$. Under this identification
we get the equalities
$$A_\nu=H^0(Y, \cO_Y(\lfloor\nu
 H\rfloor))u^\nu\quad\forall \nu\ge 0\,.$$
\esit

\bsit\label{triv1} Given a normal positively graded
$\kk$-algebra $A=\bigoplus_{\nu\ge 0} A_\nu$
a divisor $H'$ on $Y$ satisfying $A_\nu\simeq
H^0(Y, \cO_Y(\lfloor\nu H'\rfloor))$ $\forall \nu\ge 0$ can be defined
as follows (see \cite[3.5]{De}). Choose a homogeneous rational function
$u'\in (\Frac A)_1$ on $V$.
Write ${\rm div}( u'|_{V^*})= \sum_i p_i \Delta_i$, where the components
$\Delta_i$ are prime $\G_{\mathrm m} $-stable Weil divisors on  $V^*$ and
$p_i\in\Z\setminus\{0\}$. For every irreducible
component $\Delta_i$ of ${\rm div}( u'|_{V^*})$ we have
$\bar \Delta_i=\spec\, (A/I_{\bar \Delta_i})$,
where $\bar\Delta_i$ is the closure of $\Delta_i$ in $V$ and
$I_{\bar\Delta_i}$ is the graded prime ideal of $\bar \Delta_i$
in $A$. Thus the affine domain
$A/I_{\bar \Delta_i}$ is graded. Let $q_i=e(A/I_{\bar \Delta_i})$.
Then $q_i>0$ $\forall i$,
the Weil $\Q$-divisor
$H'=\sum_i \frac{p_i}{q_i}\pi_*{\Delta_i}$ satisfies
$\pi^*H'={\rm div}( u'|_{V^*})$,
and $A\simeq A(Y,H')$. Furthermore, for every component
${\Delta_i}$ the divisor
$p_i\pi_*{\Delta_i}$ is Cartier
(see \cite[Proposition 2.8]{De}).

If $A=A(Y,H)$ for a $\Q$-divisor $H$ on $Y$ then
$H'\sim H$ and $\pi^*(H'-H)={\rm div}\,(\varphi|_{V^*})$, where
$\varphi=u'/u\in ({\rm Frac}\, A)_0$. \footnote{Notice that
$({\rm Frac}\, A)_0={\rm Frac}\, A_0$ only in the case where $\dim_SY=0$.}
 \esit

\bsit\label{porog}  We keep the notation of \ref{ample}.
For some $d>0$ the $d$th
{\em Veronese subring}  of $A$,
$$A^{(d)}=A(Y,dH)=\bigoplus_{\nu\ge 0} A_{\nu d}\,,$$
is generated  over $A_0$ by  its first graded piece $A_d$:
$$A^{(d)}=A_0[A_d]$$ (see Proposition 3.3 in \cite[Ch. III, \S
1]{Bou} or \cite[Lemma 3.1.3]{Do}). This leads to an embedding of
$Y\simeq\Proj A^{(d)}$ in the projective space
$\PP_S^N\simeq\PP_{A_0}(A_d)$ and so $H$ is ample over $S$.
 Since $S$ is Noetherian, $H$ is ample over $\spec \kk$.

We call $V^{(d)}=\spec A^{(d)}$ the $d$th {\em Veronese quasicone}
associated to the pair $(Y,H)$. \esit

\bsit\label{triv2} The discussion in \S\S \ref{triv} and \ref{triv1} 
leads to the following presentations:
$$A=A(Y,H)=\bigoplus_{\nu\ge 0}
H^0(Y, \cO_Y(\lfloor\nu H\rfloor))u^\nu=\bigoplus_{\nu\ge 0}
H^0(Y, \cO_Y(\lfloor\nu H'\rfloor)){u'}^\nu =A(Y,H')\,,$$ where
$H'\sim H$ and
 $u'/u\in ({\rm Frac}\, A)_0$  (see \S \ref{triv})
 is such that $\pi^*(H'-H)={\rm div}\,(u'/u)$.
\esit
\begin{rem}[\em Polar cylinders]\label{triv3}
 We keep the
notation as in \ref{triv2}. Assume that for some non-zero
homogeneous element $f\in A_\nu$, where $\nu>0$, the open set
 $\DD_+(f)\subseteq Y$ is a cylinder i.e.,  $\DD_+(f)\simeq Z\times\A^1$
 for some variety $Z$. Then this cylinder is $H$-polar. Indeed, let $n\in\N$ be
such that $n\nu H$ is a Cartier divisor on $Y$. We have
$f^n\in A_{n\nu} = H^0(Y,\cO_Y(n\nu H))u^{n\nu}$. The rational function
$f^n u^{-n\nu}\in ({\rm Frac}\,A)_0$  being $\G_{\mathrm m} $-invariant it descends
to a rational function, say, $\psi$ on $Y$ such that
$$D:={\rm div}\,\psi+n\nu H=\pi_*{\rm div}(f^n)\ge 0\,.$$
Hence $D \in | n\nu H|$ is an effective Cartier divisor on $Y$
with $\supp D=\V_+(f^n)=\V_+(f)$. Therefore
 the cylinder
$\DD_+(f)=\spec A_{(f)}$ is $H$-polar.
\end{rem}

\eys

\bys{Generalized cones.} A quasicone  $V=\spec A$ is called a {\em generalized
cone} if $A_0=\kk$ so that $\spec A_0$ is reduced to a point.
Let us give the following example.

\begin{exa}[see e.g.\ {\cite{KPZ1}}] \label{gen-cone}
Let  $(Y,H)$ be a
polarized projective variety over $\kk$, where $H\in\Div(Y)$ is ample.
Consider the total space $\tilde V$ of the line bundle $\cO_Y
(-H)$ with zero section $Y_0 \subseteq\tilde V$. We have
$\cO_{Y_0} (Y_0) \simeq\cO_Y (-H)$ upon the natural identification
of $Y_0$ with $Y$. Hence there is a birational morphism $\varphi :
\tilde V \to V$ contracting $Y_0$. The resulting affine variety
$V=\affcone_H(Y)$ is called the {\em  generalized  affine cone over $(Y,H)$ with
vertex} $ \bar 0=\varphi(Y_0)\in V$. It comes equipped with an
effective $\G_{\mathrm m} $-action induced by the standard $\G_{\mathrm m} $-action on
the total space $\tilde V$ of the line bundle $\cO_Y (-H)$. The
coordinate ring $A=\cO(V)$ is positively graded:
$A=\bigoplus_{\nu\ge 0} A_\nu$, and the saturation index $e(A)$
equals to $1$. So the graded pieces $A_\nu$ with $\nu\gg 0$ are all
non-zero and the induced representation of $\G_{\mathrm m} $ on $A$ is
faithful. The quotient $$Y=\Proj
A=V^*/\G_{\mathrm m} \,,\quad\mbox{where}\quad V^*=V\setminus \{\bar 0\}\,,$$
can be embedded into a weighted projective space
$\PP^n(k_0,\ldots,k_n)$ by means of a system of homogeneous
generators $(a_0,\ldots,a_n)$ of $A$, where $a_i\in A_{k_i}$,
$i=0,\ldots,n$.
\end{exa}

\brems\label{gen-cone-rem} 1. Assume that $A_0=\kk$ and $V$ is normal. According to
\ref{ample}-\ref{triv2},
\be\label{grading} A=\bigoplus_{\nu\ge 0} H^0(Y, \cO_Y(\nu
H))u^\nu\quad\mbox{i.e.}\quad A_\nu=H^0(Y, \cO_Y(\nu
H))u^\nu\quad\forall\nu\ge 0\,,\ee where $u\in ({\rm Frac} A)_1$ is such that
${\rm div} (u|_{V^*})=\pi^*H$.
Since $H$ is ample   this ring is
finitely generated (see e.g.\ Propositions 3.1 and 3.2 in
\cite{Pr}).

2. If the polarization $H$ is very ample  then  $A=A_0[A_1]$ and
the affine variety $V$ coincides with the usual affine cone over
$Y$ embedded in $\PP^n$ by the linear system $|H|$. In this case
the $\G_{\mathrm m} $-action on $V^*$ is free. However, (\ref{grading}) holds
if and only if $V$ is normal that is $Y\subseteq\PP^n$ is
projectively normal. \erems

\eys

\section{The criterion}\label{section-The-criterion}
 In this section we prove our main Theorem \ref{ext-grad-int}.
Besides, in Lemma
\ref{very-ample} below we specify a range of values of
$d$ where the assertion of (b) in Theorem \ref{ext-grad-int}
can be applied.
In the sequel we fix the following setup. 

\bsit\label{hyp} Letting
$A=\bigoplus_{\nu\ge 0} A_\nu$
be a positively graded normal affine domain over $\kk$ with $e(A)=1$ we
consider the
affine quasicone $V=\spec A$ and the variety
$Y=\Proj A$ projective over the affine scheme $S=\spec A_0$.
We let  $\pi:V^*\to Y$
be the projection to the geometric quotient of $V^*$
by the natural $\G_{\mathrm m} $-action. We can write
$$A=A(Y,H)=\bigoplus_{\nu\ge 0} A_\nu,\quad\mbox{where}\quad A_\nu=
H^0(Y, \cO_Y(\lfloor\nu H\rfloor))u^\nu\,$$ with an  ample
$\Q$-divisor $H$ on $Y$ such that $\pi^*H={\rm div} (u|_{V^*})$
for some homogeneous rational function $u\in ({\rm Frac}\, A)_1$,
see \S \ref{triv2}.

Notice that in the `parabolic case' where $\dim_S Y=0$
there exists on $A$ a non-zero homogeneous
locally nilpotent derivation
`of fiber type' (that is, an $A_0$-derivation),
whatever is the affine variety $S=\spec A_0$, see
\cite[Corollary 2.8]{Li2}. In contrast, such a derivation
does not always exist if  $\dim_S Y\ge 1$.
 Hereafter we assume that $\dim_S Y\ge 1$.

Given $d>0$ we consider  the associated Veronese cone
$V^{(d)}=\spec A^{(d)}$, where $A^{(d)}=\bigoplus_{\nu\ge 0}
A_{\nu d}$. In the next example we illustrate our setting
 (without carrying the normality assumption). \esit

\bexa\label{book} In the affine space $\A^3=\spec \kk[x,y,z]$ consider the
hypersurface $$V=\V(x^2-y^3)\simeq\Gamma\times\A^1\,,$$ where
$\Gamma$ is the affine cuspidal cubic given in $\A^2=\spec \kk[x,y]$
by the same
equation $x^2-y^3=0$. Notice that $V$ is stable under the
$\G_{\mathrm m} $-action on $\A^3$ given by
$$\lambda.(x,y,z)=(\lambda^3 x,\lambda^2 y,\lambda z)\,.$$
With respect to this $\G_{\mathrm m} $-action, $\A^3$ is the generalized
affine cone over the weighted projective plane $\PP(3,2,1)$
polarized via an anticanonical divisor $H$. The divisor $H$
is ample, and $\PP(3,2,1)$ is a singular del Pezzo surface of
degree $6$. The quotient $Y=V/\G_{\mathrm m} $ is a unicuspidal   rational
curve in $\PP(3,2,1)$ with an ordinary cusp at the point
$P=(0:0:1)$. It can be polarized by an effective divisor $D$
supported at $P$ from the linear system of the restriction $H|_Y$. 
The affine surface $V\simeq \Gamma\times\A^1$ is
a cylinder, and $(Y,H)$ is cylindrical as well. The cylinder in
$Y$ consists of a single affine curve $Y\setminus \supp
D=Y\setminus \{P\}\simeq\A^1$. The natural projection $\pi:V^*\to
Y$ sends any generator $\{Q\}\times\A^1$, where $Q\in
\Gamma\setminus\{\bar 0\}$, of the cylinder $V$ onto $Y\setminus
\{P\}$.
\eexa

Recall the assertion (a) of
Theorem \ref{ext-grad-int}:

\smallskip

\noindent{\em  If a normal
affine quasicone $V=\spec A$ is cylindrical then $(Y,H)$ is}.

\smallskip

\noindent The proof given below is based on Proposition \ref{propo}.
In Corollary \ref{cyl-fac} we
start with a particular case, where the proof is rather short.
We need the following lemma.

\blem\label{spec-case} Let $R=\bigoplus_{\nu \in \Z} R_\nu$ be
a $\Z$-graded $\kk$-domain, and let $\p\in\LND(R)$ be non-zero and
homogeneous. If $e(R^{\p})$ divides $\deg \p$, then there exist
non-zero homogeneous $h \in R^{\p}$ and $f \in R^{\p}_h$ such that
$f\p \in \LND(R_{(h)})$ is non-zero.  \elem

\bproof For the proof, we use the fact that, if $t =
e(R^{\p})$ and $-d = \deg \p$, then there exists $m \in \Z$ such
that both $R^{\p}_{mt}$ and $R^{\p}_{mt+d}$ are non-zero. Up to
reversing the grading, we may suppose that $m>0$. Picking then
non-zero elements, say, $h\in R^{\p}_{mt}$ and $h_1\in
R^{\p}_{mt+d}$ we consider a homogeneous locally nilpotent
derivation  $\delta=f\p$ of degree zero on the localization $R_h$,
where $f =h_1/h \in (R^\p_h)_d$. It restricts to a locally
nilpotent derivation on $R_{(h)}$. Let us show that this
restriction is non-zero, as required. Indeed, we have\footnote{In \cite[Lemma \ref{spec-case}]{KPZ4}, the decomposition of $R_{(h)}$ was written as a direct sum decomposition. However, this is evidently false, as was pointed out to the authors by Kevin Langlois. Indeed, a nonzero element $a/h^{j}=ah/h^{j+1}$  is simultaneously contained in  $R_{tmj}h^{-j}$ and in $R_{tm(j+1)}h^{-(j+1)}$.}  
$$R_{(h)}=\sum_{j\ge 0}R_{tmj}h^{-j}\,.$$ Suppose that the restriction $\delta |_{R_{(h)}}$ is trivial. This implies that $\delta (a/h^j)= f \partial (a)/h^j=0$ for any $a \in R_{tmj}$ and $j \geqq 0$. Hence the restriction of $\partial$ to $R_{\ge 0}^{(tm)}$ is trivial, where $R_{\ge 0}^{(tm)} =\bigoplus_{j\ge 0} R_{tmj}$ is the $m$th Veronese subring of $R_{\ge 0}$. However, this is impossible
since ${\rm tr.deg} (R_{\ge 0}^{(tm)})={\rm tr.deg} (R_{\ge 0})={\rm tr.deg}
(R)={\rm tr.deg} (R^{\p})+1$.  \eproof

\bcor\label{cyl-fac} Let $A$ be a positively graded normal
affine domain over $\kk$ with $e(A)=1$, and let $\p \in \LND (A)$
be non-zero and homogeneous. Consider the presentation $A=A(Y,H)$,
where $Y=\Proj(A)$ and $H$ is an ample $\Q$-divisor on $Y$, see \S
\ref{ample}. If $e(R^{\p})\vert \deg \p$ then the pair $(Y,H)$ is
cylindrical. \ecor

\bproof Let a pair
$(A_{(h)}, f\p|A_{(h)})$ verify the conclusion of Lemma  \ref{spec-case}.
Applying to this pair the homogeneous slice construction
(see \S \ref{homog}) we obtain
a principal cylinder $\DD_+(\tilde h)=\spec A_{(h\tilde h)}$ in $\DD_+(h)$,
where $\tilde h\in \ker (f\p)\cap \im   (f\p)\subseteq A_{(h)}$ is
a  non-zero homogeneous element of degree zero.
We can write $\tilde h=a h^{-\beta}$ for some $\beta\ge 0$ and some $a\in A_\alpha$,
where $\alpha=\beta\deg h$.  Hence the cylinder
$\DD_+(\tilde h)=\DD_+(ah)= Y\setminus\V_+(ah)$ is $H$-polar,
see Remark \ref{triv3}. 
 \eproof

 The next corollary is immediate in view of Remark \ref{gen-cone-rem}.2. 

\bcor\label{addcor} Let as before $A=A(Y,H)$, and let $\p \in
\LND (A)$ be non-zero homogeneous of degree $\deg (\p)=-d$. If
$H\in\Div (Y)$ is very ample and $e (A^{\p}) =1$, then
$h^d\p\in\LND(A_{(h)})$ is non-zero for a non-zero element $h\in
A^\p_1$.  \ecor

In contrast, in case where the assumption  $e(A^\p)=1$ 
of  Corollary \ref{addcor} does not hold it is not so evident
 how one can
produce a locally nilpotent derivation on $A$ stabilizing
$A_{(h)}$ starting with a given one. Let us provide a simple
example.

\bexa\label{book1} Consider the affine plane $X=\A^2=\spec \kk[x,y]$
equipped with the $\G_{\mathrm m} $-action
$\lambda.(x,y)=(\lambda^2 x,\lambda y)$.
The homogeneous locally nilpotent derivation $\p=\frac{\p}{\p y}$ on the
algebra $A=\kk[x,y]$ graded via $\deg x=2,\,\deg y=1$ defines a principal
cylinder on $X$ with projection $x:X\to\A^1=Z$. Note that $e(A^\p)=2$.
The derivation $\p$ extends to a locally nilpotent derivation of the algebra
$$\tilde A=A[z]/(z^2-x)=\kk[z,y]\supseteq A$$ such that $e(\tilde A^\p)=1$.
The localization $A_x=\kk[x,x^{-1},y]$ extends to
$\tilde A_z=\kk[z,z^{-1},y]=\kk[z,z^{-1},s]$, where
$$s=y/z\in\tilde A_{(z)}= (\tilde A_z)_0=\kk[s]$$
is a slice of the homogeneous derivation
$\p_0=z\p\in\LND(\tilde A_{(z)})$
of degree zero. Thus $\spec \tilde A_{(z)}=\spec \kk[s]\simeq \A^1$
is a polar cylinder in $\tilde Y=\Proj \tilde A$.

The subrings $A\subseteq\tilde A$ and $A_x\subseteq\tilde A_z$
are the rings of invariants of the involution
$\tau: (z,y)\mapsto (-z,y)$ resp. $(z,s)\mapsto (-z,-s)$.
This defines the Galois $\Z/2\Z$-covers $\spec \tilde A_z\to\spec A_x$ and
$\spec \tilde A_{(z)}\to\spec A_{(x)}$.  Hence
$\spec A_{(x)}=\spec \kk[s^2]\simeq\A^1$, where $s^2=y^2/x\in A_{(x)}$,
is a polar cylinder in  $Y=\Proj A$ with a locally nilpotent derivation
$d/ds^2$. \eexa

 So in order to construct a polar cylinder for $(Y,H)$ in the general case
 one needs to apply a different strategy.  We use below the cyclic quotient
 construction (see \S \ref{fle}) 

\be\label{cycle}
Y_h=\spec F\stackrel{}{\longrightarrow}
\DD_+(h)=\spec A_{(h)}\,,\ee where $h\in A^\p_m$
is non-zero, $F=A/(h-1)$, and
$Y_h \to Y_h/\Z_m \cong \spec A_{(h)}$ is
the quotient map defined in \S 1.7. 
The key point is the following proposition.

\bprop\label{propo} Let $A=\bigoplus_{\nu\ge 0} A_\nu$
be a positively graded  normal affine domain over $\kk$,  $ Y=\Proj A$,
 and let $H$ be an ample $\Q$-divisor on $Y$ such that $A=A(Y,H)$
(see \S \ref{ample}).
 Suppose that $\dim_S Y\ge 1$, where $S=\spec A_0$.
Given a non-zero homogeneous
locally nilpotent derivation $\p\in\LND(A) $ 
there exists a homogeneous element $f\in A^\p$ such that
$\DD_+(f)=\spec A_{(f)}$ is an $H$-polar cylinder\footnote{In particular 
$\LND(A_{(f)})\neq \{0\}$. }
in $Y$.
\eprop

\bproof Let $d=-\deg\p$. We apply the homogeneous slice construction
\ref{homog}.
One can find a homogeneous element $g\in (\ker\p^2\setminus\ker\p)
\cap A_{d+m}$ such that $h=\p g\in A^\p_m$, where $m>0$ (in
particular $h$ is non-constant). Indeed, assuming to the contrary that
$A^\p\subseteq A_0$ we obtain ${\rm tr.deg} (A_0)\ge{\rm tr.deg} (A)-1$. It
follows that the morphism $Y\to S$ is finite, contrary to our
assumption that $\dim_S Y\ge 1$. Thus there exists $a\in A^\p_\alpha$,
where $\alpha>0$ and $a\neq 0$. Replacing $(g,h)$ by $(ag,ah)$, if
necessary, we may assume that $\deg h>0$. In this case the
fibers $h^*(c)$ with $c\neq 0$ are all isomorphic under the
$\G_{\mathrm m} $-action on $V=\spec A$ induced by the grading of $A$.

We use further the cyclic quotient construction, see \S \ref{fle}.
In particular, we consider the quotient \be\label{slc}
F=A/(h-1)A=A_h/(h-1)A_h= F^{\bar \p}[ \bar s],\quad\mbox{where}
\quad \bar s=g+(h-1)A\in F\ee is a slice of the induced locally
nilpotent derivation $\bar\p$ on $F$. We have $$\spec F^{\bar
\p}\simeq\V(g)\cap\V(h-1)\,,$$ where both schemes are regarded
with their reduced structure. Choosing $g$ appropriately we may
suppose that $\spec F\simeq h^*(1)$ is of positive dimension,
reduced, and irreducible. Since $\spec F\simeq\spec F^{\bar
\p}\times\A^1$ by (\ref{slc}),  then $\spec F^{\bar \p}$ is also
reduced and irreducible. Indeed, the  affine Stein
factorization (see Lemma \ref{lem: stein factorization} below) 
applied to $h$ gives $h=h_1^l$, where $m=kl$, $l\ge 1$,  and
$h_1\in A^\p_k$ is such that  the fibers $h_1^*(c)$, $c\neq 0$,
are all reduced and irreducible. Now we replace $(g,h)$ by the new
pair $(g_1,h_1)$, where $g_1=g/h_1^{l-1}\in A_h=A_{h_1}$ and
$h_1=\p g_1$. Since  the variety $\spec F^{\bar\p}$ is reduced and
irreducible $F^{\bar\p}$ is a domain. Thus we can apply Lemma
\ref{lie}.

The subgroup $\boldsymbol\mu_m\subseteq\G_{\mathrm m} $ of $m$th roots
of unity acts effectively
on $F$ stabilizing the kernel $F^{\bar\p}$. This action provides the
$\Z_m$-gradings
$$F=\bigoplus_{\sigma\in\Z_m} F_\sigma\quad\mbox{and}\quad
F^{\bar\p}=\bigoplus_{\nu=0}^{n-1} F^{\bar \p}_{[k\nu]}\,,$$ where
$m=kn$ and $k=e(A^\p)$ is such that $F^{\bar \p}_{[k\nu]}\neq 0$
$\forall \nu$, see \S \ref{fle}, Lemma \ref{lie}. The
$\boldsymbol\mu_m$-action on $F$ yields an effective
$\boldsymbol\mu_n$-action  on $F^{\bar \p}$. We have
$\bar\p:F_\sigma\to F_{\sigma-r}$, where $r=[d]\in \Z_m$.

According to (\ref{cycle})
one can write
$$A_{(h)}\simeq F^{\boldsymbol\mu_m}=F_{[0]}
=(F^{\bar \p}[\bar s])_{[0]}=\bigoplus_{j\ge 0}
 F^{\bar \p}_{[-rj]}\bar s^j\,.$$  For $\alpha\in\N$, $\alpha\gg 1$,
 one can find a non-zero element
$t\in A^\p_{e(A^\p)+\alpha m}=A^\p_{k+\alpha m}$, see Lemma
\ref{lie}. Then also $\bar t=t+(h-1)A\in F^{\bar\p}_{[k]}$ is
non-zero.
 The subgroup
 $\langle \bar t\rangle\subseteq (F_{\bar t}^{\bar \p})^\times$
 acts on $F_{\bar t}^{\bar \p}$
via multiplication permuting cyclically the graded pieces
$(F_{\bar t}^{\bar \p})_{[ik]}$, $i=0,\ldots,n-1$.
Thus $(F_{\bar t})_{[k\nu]}=(F_{\bar t})_{[0]}\bar t^\nu$ $\forall \nu$.
It follows that
$$A_{(ht)}\simeq F_{\bar t}^{\boldsymbol\mu_m}=(F_{\bar t})_{[0]}=
(F_{\bar t}^{\bar \p}[\bar s])_{[0]}=\bigoplus_{j\ge 0}
 (F_{\bar t}^{\bar \p})_{[-krj]}{\bar s}^{kj}$$$$=
 \bigoplus_{j\ge 0}(F_{\bar t}^{\bar \p})_{[0]} \left(\bar s^k\bar t^{-r}
\right)^j= \bigoplus_{j\ge 0}(F_{\bar t}^{\bar \p})_{[0]}\bar
s_1^j\,,$$ where $\bar{s_1}={\bar s}^k {\bar t}^{-r}\in
F_{[0]}$. Letting $f=ht\in A^\p_{k+(\alpha+1) m}$ we obtain
that $A_{(f)} \cong F_{({\bar t})}^{\bar \p} [{\bar s_1}]$ is a
polynomial ring. Thus $\DD_+(f)=Y\setminus \V_+(f)$ is a cylinder.
According to Remark \ref{triv3} this cylinder is $H$-polar. Now
the proof is completed.  \eproof

\noindent {\em Proof of Theorem \ref{ext-grad-int}(a)}. We
have to show that if the affine quasicone $V$ over $Y$ is
cylindrical then the pair $(Y,H)$ is. By virtue of Proposition
\ref{ext} this is true in the case where $\dim_S Y=0$, see the
discussion in \S \ref{hyp}. Otherwise the assertion follows from
Propositions  \ref{ext} and \ref{propo}. 

\medskip

This finishes the proof of part (a) of Theorem \ref{ext-grad-int}.
Part (b)  follows from the next lemma.

\blem\label{very-ample}
Assume that the pair $(Y,H)$ as in \ref{hyp}
is cylindrical with a
cylinder $$Y\setminus \supp D\simeq Z\times\A^1\,,$$
where $D$ is an effective $\Q$-divisor on $Y$ such that
$D\sim\frac{p}{q}H$ in $\Pic_\Q(Y)$ for some  coprime integers $p, q>0$.
Then the Veronese quasicone $V^{(p)}$ over $Y$
is cylindrical
and possesses a principal cylinder $\DD(h)\simeq Z'\times\A^1$, where
$Z'\simeq Z\times\A^1_*$ for some affine variety $Z$ and $h^q\in A_p$. \elem

\bproof
We have $D=\frac{p}{q}H+{\rm div} (\varphi)$ for a rational function $\varphi$ on $Y$.
Hence ${\rm div} (\varphi^q)+pH=qD\ge 0$ and so in the notation as in \ref{triv2}
$$h:=\varphi^qu^p\in A_p= H^0(Y,\cO_Y(\lfloor pH\rfloor))u^p\subseteq A\,,$$
where $u\in ({\rm Frac}\, A)_1$
satisfies ${\rm div}(u|_{V^*})=\pi^*H$. So ${\rm div}(h|_{V^*})=q\pi^*D$.
Since $$\spec (A^{(p)})_{(h)}=\DD_+(h)=Y\setminus \supp D\simeq Z\times\A^1$$ is a cylinder
we have $$(A^{(p)})_{(h)}\simeq\cO(Z)[s],\quad\mbox{where}\quad
s\in (A^{(p)})_{(h)}\quad\mbox{and}\quad  \cO(Z)\simeq (A^{(p)})_{(h)}/(s)\,.$$
Similarly as in Remark \ref{trivialization} we obtain
\be\label{noname} (A^{(p)})_h=(A^{(p)})_{(h)}[h,h^{-1}]\simeq\cO(Z)[s,h,h^{-1}]=
\cO(Z')[s]\,,\ee
where $Z'=\spec\cO(Z)[h,h^{-1}]=Z\times\A^1_*.$ Letting $\A^1=\spec \kk[s]$ we see that
$$\DD(h)=\spec\, (A^{(p)})_h\simeq Z'\times\A^1$$ is a principal cylinder in $V^{(p)}$,
as required.  \eproof

Now the proof of Theorem \ref{ext-grad-int} is completed.

\section{Final remarks and examples}\label{section-Final-remarks}
Let us start this section with the following remarks.

\brems\label{support} 1.
The assumption $D\sim \frac{p}{q}H$ of Lemma \ref{very-ample} implies the equality 
of the fractional parts 
$\{pH\}=\{qD\}$.
Hence
the  irreducible components $\Delta_i$ of the fractional part
$\{pH\}$ of the $\Q$-divisor $pH$ on $Y$ (cf.\ \ref{triv1})
are contained in $\supp \{qD\}$ and do not meet the cylinder
$Y\setminus\supp D$.

2. Suppose that  $H\in\Div Y$ is an ample Cartier divisor.
According  to Lemma \ref{very-ample} with $p=1$,
the existence of an effective divisor $D\in |H|$ such that
$Y\setminus \supp D$ is a cylinder
guarantees the cylindricity of the quasicone $V=\spec A(Y,H)$.
On the other hand, the cylindricity of $V$
does not guarantee the existence of such a divisor $D$ in
the linear system $|H|$,
but only in the  linear system $|nH|$ for some $n\in\N$, 
see Theorem \ref{ext-grad-int}(b).
We wonder whether there exists an upper  bound for such $n$ in terms
of the numerical invariants of the pair  $(Y,H)$.
This important question is non-trivial  already in the
case of del Pezzo surfaces $Y$
and pluri-anticanonical divisors $H=-mK_Y$, see
Example  \ref{exa-delPezzo} below.
\erems

Remark \ref{support}.2
together with Proposition \ref{ext}, Theorem \ref{ext-grad-int}, 
Lemma \ref{very-ample},
and  Remark \ref{curve}.2 below lead to the following corollary. 

\bcor\label{new-cor} Let $Y$ be
a normal algebraic variety over $\kk$ projective over an affine variety $S$
with $\dim_S Y\ge 1$. Let $H\in \Div (Y)$
be an ample divisor on $Y$, and
let $V=\spec A(Y,H)$ be the
associated affine quasicone over $Y$. Then
$V$ admits an effective $\G_{\mathrm a} $-action if and only
if $Y$ contains an $H$-polar cylinder. \ecor

\begin{proof} Indeed, suppose that $Y$ contains an $H$-polar cylinder.
Then by Lemma \ref{very-ample} for some $d\in\N$ the associated
Veronese quasicone $V^{(d)}$ over $Y$ admits an effective
$\G_{\mathrm a}$-action. Notice that in our setting the
$d$-sheeted cyclic cover $V \to V^{(d)}$ is non-ramified off the
vertex. Hence by Theorem 3.1 in \cite{MaMi} any effective
$\G_{\mathrm a}$-action on  $V^{(d)}$ can be lifted to $V$
 (see Remark \ref{curve}.2 below). The converse assertion follows
 immediately from Proposition \ref{ext} and Theorem \ref{ext-grad-int}.
  \end{proof}

For the details of the following examples we send the reader to
\cite{KPZ1, KPZ3}.
The latter paper inspired the present work.

\bexa\label{exa-delPezzo} The generalized cone over a smooth
del Pezzo surface $Y_d$ of degree $d$
(proper over $S=\spec \kk$)
polarized by the (integral) pluri-anticanonical divisor $-rK_{Y_d}$
admits an additive group action if $d\ge 4$
and does not admit such an action for $d=1$ and $d=2$,
whatever is $r\ge 1$. The latter follows from
the criterion of Theorem \ref{ext-grad-int}.
Indeed, in the case $d\le 2$ the pair $(Y_d,-rK_{Y_d})$
is not cylindrical (\cite{KPZ3}).
The case $d=3$ remains open.
\eexa

Remark \ref{support}.2 initiates the following definitions.

\bdefi\label{cyl-spec}
The {\em cylindricity spectrum} of a pair $(Y,H)$ is
$${\rm Sp_{cyl}} (Y,H)=\{r\in\Q_+\,|\,\exists D\in
[rH]\,\,\,\mbox{such that}\,\,\, D\ge 0\,\,\,\mbox{and}\,\,\,
Y\setminus\supp D\simeq Z\times \A^1 \}\,.$$ Clearly, ${\rm
Sp_{cyl}} (Y,H)\subseteq\Q_+$ is stable under multiplication by
positive integers. An element $r\in {\rm Sp_{cyl}} (Y,H)$ is
called {\em primitive} if it is not divisible in ${\rm Sp_{cyl}}
(Y,H)$. The set of primitive elements will be called a {\em
primitive spectrum} of $(Y,H)$. We conjecture that the primitive
spectrum is finite. \edefi

\bexas\label{exa-Platonic}
1. It may happen that  the pair $(Y,H)$ as in Theorem \ref{ext-grad-int}
is cylindrical
while the quasicone $V$ is not.
Consider, for instance, a normal generalized cone $V$
over $Y=\PP^1$, that is, a normal affine surface with a good $\G_{\mathrm m} $-action
and a quasirational singularity.\footnote{An isolated surface singularity is called
{\em quasirational} if the components of the exceptional
divisor of its minimal resolution are all rational.}
Notice that $(Y,H)$ is cylindrical
for any $\Q$-divisor $H$
on $Y$ of positive degree (see Remark \ref{0001}).
However, it was shown in \cite[Theorem 3.3]{FZ4}
that $V$ admits a $\G_{\mathrm a} $-action (that is, is cylindrical)
if and only if $V\simeq \A^2/\Z_m$ is a toric surface,
if and only if it has at most cyclic quotient singularity.
The singularities of the generalized cones
$$x^2+y^3+z^7=0
\qquad\mbox{and}\quad x^2+y^3+z^3=0$$
in $\A^3$ being non-cyclic quotient, these cones over $\PP^1$ are not
cylindrical (see \cite{FZ2}),
whereas suitable associated Veronese cones are. In terms
of the polarizing $\Q$-divisor $H$ on $Y$,
a criterion of \cite[Corollary 3.30]{Li1}  says that $V$ is cylindrical
if and only if the fractional part of $H$
is supported on at most two points of $Y=\PP^1$.
In the above examples it is supported on three points.

2. Similarly, let $a,b,c$  be a triple of positive integers coprime
in pairs, and consider the normal affine surface
$x^a+y^b+z^c=0$ in $\A^3$ with a good $\G_{\mathrm m} $-action.
According to \cite[Example 3.6]{De}
an associated $\Q$-divisor $H$ on $Y=\PP^1$ can be given as
$H=\frac{\alpha}{a}[0]+\frac{\beta}{b}[1]+\frac{\gamma}{c}[\infty]$, where
$\alpha,\beta,\gamma$ are integers satisfying
$\alpha bc+\beta ac + \gamma ab=1$.
This divisor is ample since  $\deg H=\frac{1}{abc}>0$.
For $a,b,c>1$ the fractional part of
$H$ is again supported on three points. Hence this cone,
say, $V=V_{a,b,c}$ is not cylindrical
and does not admit any $\G_{\mathrm a} $-action. At the same time
the Veronese cone $V^{(d)}$ does if and only if
at least one of the integers $a,b,c$ divides $d$.
Indeed in the latter case
the fractional part of the associated divisor
$dH$ of the Veronese cone $V^{(d)}$ is supported on at most two points.
It is easily seen that the primitive spectrum of $(\PP^1,H)$ has cardinality 3.
\eexas

\brems\label{curve} 1. Given a non-zero homogeneous derivation $\p\in\LND(A)$
of degree $d$ there exists a replica $a\p\in\LND(A^{(m)})$ of $\p$
stabilizing the $m$th Veronese subring $A^{(m)}=\bigoplus_{k\ge 0}
A_{km}$ of $A$, where $a\in A^\p_j$ for some $j\gg 0$ such that
$j+d\equiv 0\mod m$. In this way a $\G_{\mathrm a} $-action on a generalized
cone $V=\affcone_H(Y)$  induces such an action on the associated
Veronese  cone  $V^{(m)}$.
Notice that the locally nilpotent derivation on the localization
$A_h$ constructed in the proof of
Lemma \ref{very-ample} has degree zero. Hence it preserves any Veronese
subring $A_h^{(m)}$.
It follows that if $V$ is cylindrical then the associated Veronese cone
$V^{(m)}$ is cylindrical for any positive
$m\equiv 0\mod e(A^\p)$.

2. The question arises as to when a $\G_{\mathrm a} $-action on  a Veronese power
$V^{(m)}$ of a generalized cone $V=\affcone_H(Y)$  (normalized by the
standard $\G_{\mathrm m} $-action) is  induced by such an action on $V$.
The natural embedding
$A^{(m)}\hookrightarrow A$ yields an $m$-sheeted cyclic Galois cover $V\to V^{(m)}$
with the Galois group being a subgroup of the
1-torus $\G_{\mathrm m} $ acting on $V$. This cover
can be ramified in codimension 1. For instance, this is the case
if $Y$ is smooth and the ample
$\Q$-divisor $H$ is  not integral, while  $mH$ is.

In case that the cyclic cover $V\to V^{(m)}$ is unramified 
in codimension 1 the $\G_{\mathrm a} $-action on
$V^{(m)}$ can be lifted to $V$ commuting with the Galois group action
(see Theorem 1.3 in \cite{MaMi}
\footnote{See also \cite{Ka}, \cite[1.7]{FZ4}, and the proof of
Lemma 2.16 in \cite{FZ2}, where the argument must be completed.
Cf.\ Proposition 2.4 in \cite{BJ}
for a general fact on lifting algebraic group actions to
an \'etale cover over a complete base
in arbitrary characteristic.}).
\erems

The following simple example\footnote{Cf.\ \cite[Example
2.17]{FZ2}. The double definition of $A$ in  \cite[Example
2.17]{FZ2} is not correct; the correct definition is given by the
second equality $A=\bigoplus_{\nu\neq 1} A_\nu$, while for our
purposes the first equality is more suitable.} shows that without the normality
assumption for the quasicone $V$, it is impossible in general to
lift  to $V$ a given $\G_{\mathrm a} $-action on a Veronese cone  $V^{(m)}$.

\bexa\label{exa-non-normal} Consider the polynomial algebra
$\tilde A=\kk[x,\,y]$ with the standard grading and a homogeneous locally nilpotent
derivation $\partial=y {\p \over\p x}$ of degree $0$.
Consider also
a non-normal subring
$$\tilde B=\kk[x^2,\,xy,\,y^2,\,x^3,\,y^3]\subseteq \tilde A\,$$
with normalization $\tilde A$. Note that $\p$ does not stabilize
$\tilde B$, 
while $y^3 {\p \over\p x}$ does stabilize $\tilde{B}$,
i.e., $\tilde{B}$ is not a rigid
ring.

On the other hand, the involution $\tau:(x,y)\mapsto
(-x,-y)$ acts on $\tilde A$ leaving $\tilde B$ invariant.
Furthermore, letting $G=\langle\tau\rangle\simeq\Z/2\Z$ we obtain
$$\tilde A^G=\tilde B^G=\bigoplus_{\nu=0}^\infty
\tilde A_{2\nu}=:A\,.$$ Let $\tilde V=\spec \tilde B$ and $V=\spec
A$. Then $\tilde V=\affcone (\tilde \Gamma)$ is  a generalized
affine cone over the smooth projective rational  curve $\tilde
\Gamma\subseteq\PP^4(2,2,2,3,3)$ given by $(x:y)\mapsto
(x^2:xy:y^2:x^3:y^3)$, while $V=\affcone (\Gamma)$ is  the usual
quadric cone over a smooth conic $\Gamma\subseteq\PP^2$.
The embedding $A\hookrightarrow \tilde B$ induces a 2-sheeted
Galois cover $\tilde V\to V$ ramified only over the vertex of $V$.
The derivation $\p$ stabilizes $A$, and the induced
$\G_{\mathrm a} $-action
on $V$ lifts to the normalization $\A^2$ of $\tilde V$, and also
to $\tilde V^*=\tilde V\setminus\{\bar 0\}\simeq\A^2\setminus
\{\bar 0\}$. However, since $\p$ does not stabilize $\tilde B$
this action cannot be lifted to the cone $\tilde V$. \eexa

\brem[\em Affine Stein factorization]\label{Stein
factorization}  In the proof of Proposition \ref{propo} we have
used the following affine version of the classical Stein
factorization. It should be well known; for the lack of a
reference we provide a short argument. \footnote{We are
grateful to Hubert Flenner for indicating this argument. An
alternative proof, also discussed with him, consists to define
$Y'$ below as spectrum of the normalization of the algebra
$\cO(Y)$ in $\cO(X)$. In this way we avoid the desingularization,
but the proof becomes somewhat longer. } \erem

\blem\label{lem: stein factorization} Given a dominant morphism
$f: X \to Y$ of affine varieties there exists a decomposition
$f=g\circ f'$, where $f': X \to Y'$ is a morphism with irreducible
general fibers and the morphism $g: Y' \to Y$ is
finite\footnote{So $Y'$ is affine. }.\elem

\bproof Compactifying $X$ and $Y$ and resolving indeterminacies of
the resulting rational map we can extend $f$ to a morphism of
projective varieties $\tilde f:\tilde X \to \tilde Y$. We then
restrict $\tilde f$ over $Y$  to get a proper morphism $\bar
f:\bar X\to Y$, where $\bar X={\tilde f}^{-1}(Y)$ is open in
$\tilde X$. Now by \cite[Ch.\ III, Cor.\ 11.5]{Ha} there is a
factorization $\bar f=g\circ \bar f'$, where $\bar f':\bar X\to
Y'$ is a morphism with irreducible general fibers and the morphism
$g: Y' \to Y$ is finite. Letting $f'=\bar f'|_X$ we are done. 
\eproof

\end{document}